\documentclass{amsart}
\usepackage[latin1]{inputenc}
\usepackage{amssymb}
\usepackage{amsfonts}
\usepackage{amsmath}
\usepackage{amsthm}
\usepackage{latexsym}
\usepackage{bm}
\usepackage{bbm}
\usepackage{mathrsfs}
\usepackage[T1]{fontenc}
\usepackage{ae}
\usepackage[english]{varioref}
\usepackage{enumitem}

\newtheorem{thm}{Theorem}[section]
\newtheorem*{thm*}{Theorem} 
\newtheorem{lemma}[thm]{Lemma}

\theoremstyle{definition}
\newtheorem{defn}[thm]{Definition}

\theoremstyle{remark}
\newtheorem{remark}[thm]{Remark}

\newcommand{\N}{\mathbbm{N}}
\newcommand{\Z}{\mathbbm{Z}}

\newcommand{\R}{\mathbbm{R}}
\newcommand{\C}{\mathbbm{C}}
\newcommand{\Hy}{\mathbbm{H}}
\newcommand{\een}{\mathbbm{1}}

\newcommand{\dr}{\textup{SL}(d,\R)}
\newcommand{\dz}{\textup{SL}(d,\Z)}
\newcommand{\iir}{\textup{SL}(2,\R)}
\newcommand{\iiz}{\textup{SL}(2,\Z)}
\newcommand{\drdz}{\dr/\dz}
\newcommand{\dd}{\textup{\textbf{d}}}

\newcommand{\para}[1]{\left(#1\right)}
\newcommand{\kpara}[1]{\left[#1\right]}
\newcommand{\num}[1]{\left \vert #1\right \vert}
\newcommand{\tub}[1]{\left \{#1\right \}}
\newcommand{\norm}[1]{\left\Vert #1\right\Vert}
\newcommand{\skr}{\mathfrak}
\newcommand{\tyk}{\mathcal}
\newcommand{\supp}{\textup{supp}}
\newcommand{\Ad}{\textup{Ad}}

\renewcommand{\phi}{\varphi}
\renewcommand{\epsilon}{\varepsilon}

\numberwithin{equation}{section}

\begin{document}

\title[EVD's for one-parameter actions on homogeneous spaces]{Extreme value distributions for one-parameter actions on homogeneous spaces}

%    Remove any unused author tags.

%    author one information
\author{Maxim S\o lund Kirsebom}
\address{Department of Mathematics, University of Hamburg, Bundesstrasse 55, 20146 Hamburg, Germany}
\curraddr{}
\email{maxim.kirsebom@uni-hamburg.de}
\thanks{This research was supported by ERC grant 239606 and DFF/Marie Curie COFUND grant 5051-00184}



\date{02.04.2015}

\dedicatory{}

\begin{abstract}
In this paper we study extreme value distributions for one parameter actions on homogeneous spaces of Lie groups. We study both shortest vectors in unimodular lattices, maximal distance excursions and closest distance returns of a one-parameter action. For certain sparse subsequences of the one-parameter action and by taking the maximum over a moving interval of indices we prove non-trivial estimates for the limiting distribution in all cases.
\end{abstract}

\maketitle

\section{Introduction}

The setting of this paper is a one-parameter group acting on a non-compact, finite volume homogeneous space of a Lie group. We evaluate an observable along orbits of the system and consider the maximum and minimum values obtained up to a given time. The aim is to understand the distribution of this maximum and minimum as time goes to infinity.

Before elaborating on the details we give the broader context into which this work fits.

\subsection{Shrinking targets, logarithm laws and Dynamical Borel-Cantelli Lemmas}

Much attention has recently been devoted to so-called shrinking target problems. Heuristically, these may be described as questions probing the recurrence properties of dynamical systems to sequences of sets that shrink in some given sense as time increases. The term was first coined by Hill and Velani in \cite{HiVe}, but the interest in this class of problems goes further back.

A popular question to ask in this direction is whether a system obeys a so-called logarithm law. To the authors knowledge, the first such law was formulated and proven by Sullivan in \cite{Sull}. His theorem concerns maximal excursions into a cusp by the geodesic flow on certain hyperbolic manifolds.
\begin{thm*}[Sullivan's Logarithm Law, \cite{Sull} Theorem 2]%\label{Sullivan}
	Let $V=\Hy^{d+1}/\Gamma$ where $\Gamma$ is a discrete subgroup of hyperbolic isometries such that $V$ has finite volume but is non-compact. Let $\textup{dist}\,v(t)$ denote the distance from a fixed point in $V$ to the point achieved after traveling a time $t$ along the random geodesic with initial direction $v$. Then for almost all starting directions $v$ of geodesics
	\begin{equation*}
		\limsup_{t\to\infty} \frac{\textup{dist}\,v(t)}{\log t}=\frac{1}{d}.
	\end{equation*}   
	\end{thm*}
In this context, shrinking targets are understood as shrinking neighborhoods of the cusps (i.e. complements of balls of increasing radius). In words, Sullivan's Logarithm Law tells us that almost surely, the maximal distance the geodesic flow has ventured away from a fixed point at time $t$ is asymptotically of the same order as $\log t$. Logarithm laws turn out to hold for a vast variety of dynamical systems, see for example \cite{AGP1}, \cite{AtMa1}, \cite{Gala}, \cite{HePa}, \cite{KleinMarg}, \cite{Mark}, \cite{Mauc}, \cite{StVe}, \cite{Tseng}, for a (non-exhaustive) list of examples in different settings. 

While a logarithm law provides interesting information about a system's recurrence to shrinking targets, it is natural to ask whether we can say something more precise. Very often this is the case, and indeed, many logarithm laws are not proven directly, but appear as corollaries of more general statements known as dynamical Borel-Cantelli Lemmas. Let $\tub{A_n}$ denote a sequence of sets in a probability space $(X,\mu)$. The classical Borel-Cantelli Lemma then states that

\begin{enumerate}
\itemsep1em
\item $\sum \mu(A_n)<\infty$ implies $\mu\para{\tub{x:x\in A_n \text{ for }\infty \text{ many }n}}=0$.
\item Assume that the sets $A_n$ are independent. Then $\sum \mu(A_n)=\infty$ implies $\mu\para{\tub{x:x\in A_n \text{ for }\infty \text{ many }n}}=1$.
\end{enumerate}

Given a dynamical system $(X,T,\mu)$ we ask if almost surely $x\in T^{-n} A_n$ (i.e. $T^n x\in A_n$) for infinitely many $n$. Often, though not always (see \cite{Fayad}), we may answer this in the affirmative by replacing the assumption of independence with some mixing property of the dynamical system. When this is the case we say that we have a dynamical Borel-Cantelli Lemma for the sequence $\tub{A_n}$. A logarithm law then typically follows if we can prove a dynamical Borel-Cantelli Lemma for a sufficiently large family of sequences $\tub{A_n}$ satisfying $\sum \mu(A_n)=\infty$.

\subsection{Extreme Value Laws in dynamics}

The focus of this paper however, is on a different generalization of a logarithm law known as an extreme value law (EVL) or extreme value distribution (EVD). In this context it may be viewed as a distributional sibling to the "almost-sure" approach to shrinking target problems. The theory behind this is known as extreme value theory (EVT) and attempts to answer questions of the following kind. Let again $(X,T,\mu)$ denote a dynamical system and $\tyk{D}:X\to\R$ a measurable function which we will also refer to as an observable. In the setting of Sullivan, $\tyk{D}$ could for example be chosen as the distance to a fixed point, i.e. $\tyk{D}(\cdot)=\textup{dist}(\cdot,x_0)$, $x_0\in V$. Define $M_n:X\to\R$ by
\begin{equation*}
	M_n(x):=\max_{0\leq i \leq n} \tyk{D}(T^i x).
\end{equation*}
We are interested in the distribution of this maximum (under suitable normalization) when $n$ becomes large. That is, we ask if there exists real sequences $a_n>0$ and $b_n$ such that
\begin{equation}\label{EVLconv}
	\mu\para{\tub{x\in X: \frac{M_n(x)-b_n}{a_n}\leq r}}\to G(r) \enskip\text{for}\enskip n\to\infty
\end{equation}
for some non-trivial distribution function $G$. We say that $G$ is trivial if it only takes the values 0 and 1. If the non-trivial limit exists we are interested in determining the form of the limit as well. Clearly the same question may be asked for a continuous flow on $X$ with the definition of the maximum adapted accordingly.

If the $\tyk{D}(T^i x)$ were independent as random variables, the answers to the above questions would be well understood. Indeed, the Extremal Types Theorem (see \cite{LLR}) says that for i.i.d. sequences, $G$ can only attain one of three forms. More precisely, if there exist sequences $a_n>0$ and $b_n$ such that \eqref{EVLconv} holds for some non-trivial $G$, then $G$ is on one of the following three forms up to type (we say that $G_1$ and $G_2$ are of the same type if $G_1(r)=G_2(ar+b)$ for some $a>0,b\in\R$):
	\begin{align*}
	&\text{\textbf{Type 1 (Gumbel):}}\quad G(r)=\exp(-\exp(-r))\;,\;r\in\R\\
	&\text{\textbf{Type 2 (Fr\'echet):}}\quad G(r)=
	\begin{cases}
	0 & \text{for } r\leq 0\\
	\exp(-r^{-\alpha}) & \text{for } r>0 
	\end{cases}
	\quad \alpha>0
	\\
	&\text{\textbf{Type 3 (Weibull):}}\quad G(r)=
	\begin{cases}
	\exp(-(-r)^{\alpha}) & \text{for } r\leq 0\\
	1 & \text{for } r>0
	\end{cases}
	\quad \alpha>0.
	\end{align*}
	
Furthermore, we also have a necessary and sufficient condition for convergence to a non-trivial limit (again, see \cite{LLR}). Let $\tau:\R\to\R$ and let $u_n(r)$ be a sequence of real functions. Then
\begin{align}
\label{EVT1}
n\mu\para{\tub{x:\tyk{D}(x)>u_n(r)}}\to\tau(r) \iff \mu(\tub{x:M_n\leq u_n(r)})\to e^{-\tau(r)}
\end{align}
as $n\to \infty$. We refer to $\mu\para{\tub{x:\tyk{D}(x)>z}}$ as the tail-distribution function of $\tyk{D}$.

It is well known that dynamical systems do not typically generate independent random variables. If the system is mixing, this may be understood as an asymptotic independence and, as was the case for the dynamical Borel-Cantelli Lemma, we may hope to obtain similar results by replacing independence with mixing. The general framework of EVT in the presence of dependence is described in \cite{LLR}, Chapter 3. In the past two decades great progress has been achieved in applying this framework to certain dynamical systems with various hyperbolicity assumptions. It began with the paper by Collet \cite{Collet} who proved an EVL for closest returns to a fixed point for a non-uniformly hyperbolic $C^2$ interval map. Implicitly in his proof, Collet improved the general framework for EVT in the presence of dependence. His idea was formalized and developed by Freitas and Freitas in \cite{FF} and subsequently applied in many different cases by several authors. See \cite{Freitas} and references therein for a survey of developments in this field between Collet's paper and 2013.

\subsection{Homogeneous dynamics}

While the technical details are not yet important, we might as well let the reader in on the general setting of the paper. Definitions and explanations will follow in later sections.

Let $G$ denote a connected semi-simple Lie group with finite center and no compact factors. Let $\Gamma<G$ be an irreducible lattice such that $G/\Gamma$ is not compact. Set $X=G/\Gamma$ and let $a_t:=\tub{a_t : t\in \R}$ denote a partially hyperbolic one-parameter subgroup of $G$. Let $\mu$ denote the Haar measure on $X$ and notice that $\mu(X)<\infty$. This allows us to further assume that $\mu$ is normalized. i.e. $\mu(X)=1$.

This is essentially the setting of \cite{KleinMarg} which represents one of the first treatments of shrinking target problems in homogeneous dynamics and served as the inspiration for this work. In \cite{KleinMarg}, Kleinbock and Margulis proved dynamical Borel-Cantelli Lemmas and established a generalization of Sullivan's Logarithm Law to homogeneous spaces. Since then much progress has been achieved in this direction, see for example \cite{Kelm}, \cite{KeYu}, \cite{KleinZhao} for some interesting recent works. See also \cite{Athr} for a well written survey on shrinking target problems and logarithm laws in the setting of homogeneous dynamics.

Unfortunately, very little is yet known about EVT for homogeneous dynamical systems. In particular in the non-compact finite-volume case, where the only known result appears to be due to Pollicott \cite{MarkP} in the special case of the geodesic flow on $\iir/\iiz\simeq T^1(\Hy^2/\iiz)$.
\begin{thm*}[Pollicott, \cite{MarkP} Theorem 2]
	Let $V=\Hy^2/\iiz$, let $\mu_L$ denote the Liouville measure on $T^1 V$, let $h$ denote the function returning the hyperbolic height above the horizontal line $\textup{Im}(z)=1$ and let $\gamma_t$ denote the geodesic flow on $T^1 V$. For any $r\in\R$ we have
	\begin{align*}
	\lim_{T\to\infty} \mu_L\para{(z,v)\in T^1V:\max_{0\leq t\leq T} h(\gamma_t(z,v)) \leq r+\log T}=e^{-\frac{3}{\pi^2}e^{-r}}.
	\end{align*} 
\end{thm*}
\begin{remark}
	The theorem in \cite{MarkP} is stated with $6/\pi$ instead of $3/\pi^2$. However, this was a typo. 
	\end{remark}
Pollicott also showed how a logarithm law follows as a corollary of the above theorem, justifying the notion that EVL's are generalizations of logarithm laws. Indeed, as we demonstrated in \cite{Maxim1}, even an upper bound on the $\limsup$ for the distribution of maxima along a sparse subsequence is sufficient to recover a logarithm law.
Pollicott's proof of the theorem uses an EVL for the digits of continued fractions expansions of real numbers (see \cite{Galambos}) as well as the connections between geodesics on $\Hy^2$ and continued fractions. These ingredients are not known in higher dimensions and hence the proof can not easily be adapted to higher generality.
The approach of Collet to proving EVL's in dynamical systems is also not applicable in our setting. The very short explanation is that Collet's approach relies heavily on the ability to model the dynamical system in question by Young Towers. This is not possible in the setting of our interest and hence this approach is also not of help.

We finish this introduction with a brief description of the strategy employed in this paper. Without proven techniques to build on we take a direct approach to the problem. The assumption that $a_t$ is partially hyperbolic ensures that our system is exponentially mixing with respect to Sobolev observables. This, first and foremost, makes it reasonable to suspect from a statistical point of view, that the system behaves similar to if it were independent. For simplicity, say that we want to establish an EVL for the time-one map of $a_t$. Let $\tyk{D}:X\to \R$ denote an observable, and write $u_n=a_n r+b_n$ for some sequences $a_n>0$ and $b_n$. The starting point is the following trivial, but important rewriting of the term we want to estimate. 
\begin{align*}
\mu\para{x\in X:\max_{1\leq i\leq n}\tyk{D}(a_i x)\leq u_{n}}&=\mu\para{\bigcap_{i=1}^n\tub{x\in X:\tyk{D}(a_i x)\leq u_{n}}}\\
&=\int_X \een_{\bigcap_{i=1}^n\tub{x\in X:\tyk{D}(a_i x)\leq u_{n}}}d\mu(x) \\
&=\int_X \prod_{i=1}^n\een_{\tub{y\in X:\tyk{D}(y)\leq u_{n}}} (a_i x)d\mu(x).
\end{align*}
The idea is then to approximate the characteristic functions by Sobolev functions and apply exponential mixing $n-1$ times to obtain a product of $n-1$ integrals plus the sum of $n-1$ error terms. As it turns out, the exponential rate of decay is not sufficient to neutralize the $n-1$ error terms as $n$ goes to infinity. To obtain this we are forced to consider sparse subsequences of times (though not necessarily integer) as well as maxima over a moving window of indices. This is the simple idea that underpins this article. Details will be clear in the next section where main results are stated.

\subsection*{Acknowledgements} The author would like to thank Alex Gorodnik for many helpful discussions and comments concerning this paper.

\section{Main results}
\label{submainres}

We assume the general setting described in the introduction and introduce the rest of the notation necessary for stating the main theorems. Let $\tub{m_j}_{j=0}^{\infty}\subset\R_{\geq 0}$ denote a strictly increasing sequence and let $\tub{\alpha_n},\tub{\beta_n}$ be sequences of natural numbers with $\alpha_n< \beta_n$ for all $n\in\N$. Set $I_n=\tub{m_{\alpha_n}, m_{\alpha_n+1},\dots,m_{\beta_n}}$. $\tyk{D}$ will be used to denote a general observable on $X$. We define the maximum over times in $I_n$ by
\begin{align}\label{maxIn}
M_{I_n}(x):=\max_{i\in I_n} \tyk{D}(a_i x)=\max_{\alpha_n\leq j \leq \beta_n} \tyk{D}(a_{m_j}).
\end{align}
This is the maximum which we will consider throughout the paper while varying the choice of observable $\tyk{D}$ and the conditions on $\alpha_n, \beta_n$ and $m_j$.  
Generally we require that $\alpha_n\to\infty$ as well as $\beta_n-\alpha_n\to\infty$. In each of the main theorems we make further requirements on the growth rate of $m_j$ and $\beta_n-\alpha_n$ depending on the setting. We also fix the notation
\begin{align*}
N_n:=\beta_n-\alpha_n+1.
\end{align*}
throughout.

Henceforth we will denote $\limsup$ and $\liminf$ by $\varlimsup$ and $\varliminf$. We will also use shorthand notation for measures of sets, for example writing
\begin{align*}
	\mu(M_n\leq u_n(r)):=\mu\para{\tub{x\in X:M_n(x)\leq u_n(r)}}.
\end{align*}
Let us also here mention that we will often write $u_n:=u_n(r)$ when considering scaling sequences and the dependence on $r$ is not important.

We shall make use of big-$O$ and small-$o$ notation in the standard way and use $O_x$ to mean that the implicit constant depends on $x$. When it suits the situation we will also use the Vinogradov symbol $\ll$ instead of the big-$O$ notation and $\ll_x$ will mean that the implicit constant depends on $x$.

Generally, all integrals in the paper will be taken over $X$ with respect to $\mu$ unless explicitly stated otherwise. We use this to simplify notation by writing 
\begin{equation*}
\int f:=\int_X f d\mu.
\end{equation*}
This should not lead to any confusion and will shorten some rather long equations later in the paper.

\subsection{Shortest vectors on the space of unimodular lattices}
\label{subsecmain2}
We first consider a special case of our general setting. Let $X=\tyk{L}_d$ denote the space of unimodular lattices in $\R^d$ and $\mu$ the normalized Haar measure on $\tyk{L}_d$. Recall that $\tyk{L}_d$ can be identified with $\drdz$. Let $a_t:=\tub{a_t : t\in\R}$ denote a one-parameter subgroup of $\dr$ and assume that $a_t$ is partially hyperbolic. Let $\Delta:\tyk{L}_d\to\R$ be given by
\begin{align}
\label{Deltadef}
\Delta(\Lambda)=\max_{v\in \Lambda\backslash\tub{0}}\log\para{\frac{1}{\norm{v}}}.
\end{align}
The observable $\Delta$ plays an important role in the connections between flows on $\tyk{L}^d$ and Diophantine approximation. Note that up to a change of variables $\Delta$ returns the length of the shortest non-zero vector in the lattice. Set $\tyk{D}:=\Delta$ and $u_n=r+\frac1d\log N_n$. We obtain the following result.
\begin{thm}[]
	\label{thm:[ExLattice]}
	There exists an explicit constant $C\in (0,1]$ (see Theorem \ref{thm:[GenSDL]}) such that by assuming
	\begin{align}
	\label{thm43cond}
	\sup_{j\in\N}\frac{m_{j-1}}{m_j}<C
	\end{align}
	we get
	\begin{align*}
	\lim_{n\to\infty} \mu\para{M_{I_n}\leq u_n}=e^{-\frac{V_d}{2\zeta(d)}e^{-dr}},
	\end{align*}
	where $V_d$ denotes the volume of the unit ball in $\R^d$.
\end{thm}

\subsection{Maximal distance and closest returns on homogeneous spaces of Lie groups}
\label{subsecmain1}

We restate the general setting.
Let $G$ denote a connected semi-simple Lie group with finite center and no compact factors. Let $\Gamma<G$ be an irreducible lattice such that $G/\Gamma$ is not compact. Set $X=G/\Gamma$ and let $\mu$ denote the normalized Haar measure on $X$. Let $\dd$ be the Riemannian metric on $X$ chosen by fixing a right invariant Riemannian metric on $G$ which is bi-invariant with respect to a maximal compact subgroup of $G$. Let $d$ denote the dimension of $G$. 

Let $a_t:=\tub{a_t : t\in \R}$ denote a one-parameter subgroup of $G$ and assume that $a_t$ is partially hyperbolic.
Let $u_n=r+\frac1v\log N_n$ and for a fixed point $x_0\in X$ set $\tyk{D}(\cdot):=\dd(\cdot,x_0)$. 

We prove the following for maximal distance excursions.
\begin{thm}[]
\label{thm:[exRieDist]}
There exists an explicit constant $C\in (0,1]$ (see Theorem \ref{thm:[GenSDL]}) and strictly positive constants $w_1,w_2$ and $v$, such that by assuming
\begin{align}
\label{thm41cond}
\sup_{j\in\N}\frac{m_{j-1}}{m_j}<C,
\end{align}
we get that for every $x_0\in X$
\begin{align*}
e^{-w_1e^{-vr}}\leq \varliminf_{n\to\infty} \mu&\para{M_{I_n}\leq u_n}\leq\varlimsup_{n\to\infty} \mu\para{M_{I_n}\leq u_n}\leq e^{-w_2e^{-vr}}.
\end{align*}
\end{thm}

Set now $\tyk{D}(\cdot):=-\log\dd(\cdot,x_0)$ for some fixed point $x_0$ and notice that $\tyk{D}$ becomes large when $\dd(\cdot,x_0)$ becomes small. Set also $u_n=r+\frac1d\log N_n$.
We prove the following for closest distance returns.
\begin{thm}[]
\label{riedistret}
There exists explicit constants $C\in (0,1]$ and $\sigma>0$ (see Theorem \ref{thm:[ClosRet]}) and a strictly positive constant $\kappa$, such that by assuming
\begin{align}
\label{thm42cond}
\frac{1}{\rho}=:\sup_{j\in\N} \frac{m_{j-1}}{m_j}<C
\end{align} 
and
\begin{align*}
N_n=o\para{e^{\sigma \rho^{\alpha_n}}},
\end{align*}
we get that for every $x_0\in X$
\begin{align*}
\lim_{n\to\infty} \mu\para{M_{I_n}\leq u_n}=e^{-\kappa e^{-d r}}.
\end{align*} 
\end{thm}

A simple example for which the conditions of Theorem \ref{riedistret} are satisfied is $\beta_n=2n$, $\alpha_n=n$ and $m_j=q^j$ for some sufficiently large $q\in\R$. For some $\kappa>0$ we then have for all $x_0\in X$
\begin{align*}
\lim_{n\to\infty} \mu\para{\max_{n\leq j\leq 2n} -\log\dd(a_{q^j}x,x_0)\leq u_n}=e^{-\kappa e^{-d r}}.
\end{align*}

\subsection{Structure of the paper}
We begin in Section \ref{Prem} by describing the main tool used in the proofs, namely exponential mixing. Thereafter the paper splits in two parts due to the following consideration. The main results listed above concern the three observables
\begin{enumerate}[label=(\roman*)]
	\item\label{D1} $\Delta(\Lambda)=\max_{v\in \Lambda\backslash\tub{0}}\log\para{\frac{1}{\norm{v}}}$, on $X=\drdz$.
	\item\label{D2} $\dd(x,x_0)$, on $X=G/\Gamma$.
	\item\label{D3} $-\log \dd(x,x_0)$, on $X=G/\Gamma$.
\end{enumerate}
Case \ref{D1} and \ref{D2} can be treated almost identically while case \ref{D3} is treated differently. The main distinguishing property of the observables is that \ref{D1} and \ref{D2} are uniformly continuous while \ref{D3} is not.

For this reason, Section \ref{UniContObs} is dedicated to the case of uniformly continuous observables. In this section we introduce necessary tools and prove a general theorem from which Theorem \ref{thm:[ExLattice]} and Theorem \ref{thm:[exRieDist]} will follow.
In Section \ref{ClosRet} we prove Theorem \ref{riedistret}. Many of the tools and preliminary results necessary for the proof of Theorem \ref{riedistret} are developed in Section \ref{UniContObs}. Hence Section \ref{ClosRet} mainly highlights the differences arising from the lack of uniform continuity of the observable.  

\section{Exponential mixing}
\label{Prem}
We henceforth assume the setting of Subsection \ref{subsecmain1}.
The assumption that $a_t$ is partially hyperbolic means that $\Ad(a_1)$ has at least one eigenvalue different from 1 in absolute value. Since $G$ is semi-simple it follows in particular that at least one eigenvalue of $\Ad(a_1)$ is strictly bigger than 1 in absolute value. 
Denote the eigenvalues of $\Ad(a_1)$ by $\lambda_1,\dots,\lambda_d\in \C$ and for later convenience we write $\num{\lambda_j}:=e^{\gamma_j}$ for some $\gamma_j>0$, for $1\leq j\leq d$ and set 
\begin{equation}
\label{gamma}
\gamma:=\max_{1\leq j\leq d}\gamma_j.
\end{equation}
Then the largest eigenvalue in absolute value is $e^{\gamma}=\max_{1\leq j\leq d} \num{\lambda_j}$.
Fix a Jordan basis $B=\tub{\zeta_1,\dots,\zeta_d}$ for $\Ad(a_1)$ such that the matrix representing $\Ad(a_1)$ is a real Jordan form.

Decay of correlations plays a central role in developing statistical limit theorems for dynamical systems.
The type of mixing known in our setting is of exponential rate for Sobolev observables. To precisely state this mixing property we need to define Sobolev norms. First we define what we mean by the derivative of a smooth function on $X$.

\begin{defn}
	The derivative of a function $f\in C^{\infty}(X)$ in the direction of an element $\zeta\in \skr{g}$ will be denoted $D_{\zeta}f$ and is defined by
	\begin{align*}
	(D_{\zeta}f)(x)=\frac{d}{dt} f(\exp{(t\zeta)}x)|_{t=0}.
	\end{align*}
\end{defn}
\begin{remark}
	We define the derivative of a function $f\in C^{\infty}(G)$ in the direction of $\zeta\in\skr{g}$ analogously.
\end{remark}
Note that we may think of $D_{\zeta}$ as a differential operator on $C^{\infty}(X)$. It satisfies $D_{a\zeta_1+b\zeta_2}=aD_{\zeta_1}+bD_{\zeta_2}$ for all $a,b\in\R$ and any $\zeta_1,\zeta_2\in\skr{g}$. Also, the differential operators commute. Hence any set of the form $W=\tub{(\zeta_1,n_1),
	\dots,(\zeta_d,n_d)}$, $n_i\in\N$ gives rise to a higher order differential operator through composition, i.e.
\begin{align*}
D_W=D_{\zeta_1}^{n_1}\cdots D_{\zeta_d}^{n_d}.
\end{align*}
We define the degree of $D_W$ as $\textup{deg}(D_W):=n_1+\dots+n_d$.  
\begin{defn}[Sobolev norm]
	For an integer $k\geq 1$ and $f\in C^{\infty}(X)$, the "$L^2$, degree $k$" Sobolev norm of $f$ is denoted $S_k(f)$ and is given by
	\begin{align*}
	S_k(f)^2=\sum_{\textup{deg}(D_W)\leq k}\norm{D_W f}_2^2.
	\end{align*}
\end{defn}
\begin{remark}
	Clearly, the definition of the Sobolev norm depends on the choice of basis for $\skr{g}$, but a change of basis only changes $S_k(f)$ by a bounded factor. We will only consider functions with finite Sobolev norm so for our purposes it is safe to omit this dependence in the definition.
\end{remark}
Denote by $\tyk{S}_k(X):=\tub{f\in C^{\infty}(X):S_k(f)<\infty}$. The exponential mixing property on homogeneous spaces has a long history.
In this work, we use the following version.
\begin{thm}[\cite{KleinMarg}]
	There exist constants $\delta>0$, $C>0$ and $k\in\N$ such that for any two functions $f,g\in \tyk{S}_k(X)$ and for any $t\geq 0$ we have
	\begin{align}
	\label{expomix2}
	\num{\int f(x)g(a_t x)-\int f \int g } \leq Ce^{-\delta t} S_k(f) S_k(g).
	\end{align}
\end{thm}
The proof of this estimate appears in \cite{KleinMarg} under the condition of having a strong spectral gap which was later verified in \cite{KeSa}.

\section{EVD's for uniformly continuous observables}
\label{UniContObs}

Assume throughout this section that $\tyk{D}:X\to\R$ is uniformly continuous. For the purpose of this work, the most important property of $\tyk{D}$, along with the uniform continuity, is how well we know the asympotics of its tail distribution function. This motivates the following definition.
\begin{defn}[]
	\label{DistanceLike}
	\textbf{(DL)} For strictly positive constants $w_1,w_2$ and $v$, we say that $\tyk{D}$ is $(w_1,w_2,v)$-DL ("Distance-Like") if it is uniformly continuous and satisfies
	\begin{align*}
	w_1 e^{-vz}\leq \mu\para{x:\tyk{D}(x)\geq z}\leq w_2 e^{-vz}\enskip,\enskip\forall\, z\in \R.
	\end{align*}
	\textbf{(SDL)} For strictly positive constants $w$ and $v$, we say that $\tyk{D}$ is $(w,v)$-SDL ("Strong-Distance-Like") if it is uniformly continuous and satisfies
	\begin{align*}
	\mu\para{x:\tyk{D}(x)\geq z}=w e^{-vz}+o(e^{-vz}) \quad\text{ as } z\to\infty.
	\end{align*}
\end{defn}
The notion of distance-like functions was introduced by Kleinbock and Margulis in \cite{KleinMarg}. In the same paper they proved that $\dd(\cdot,x_0)$ is DL in the setup of Theorem \ref{thm:[exRieDist]} (\cite{KleinMarg} Proposition 5.1) and that $\Delta$ is SDL in the setup of Theorem \ref{thm:[ExLattice]} (\cite{KleinMarg} Proposition 7.1). With this in mind it is clear why the following theorem immediately implies Theorem \ref{thm:[ExLattice]} and Theorem \ref{thm:[exRieDist]}.
\begin{thm}[]
	\label{thm:[GenSDL]}
	Assume that $\tub{m_j}$ satisfies
	\begin{align*}
	\sup_{j\in\N}\frac{m_{j-1}}{m_j}<\min\para{1,\frac{\delta}{k\gamma'}},
	\end{align*}
	where $k,\delta$ are given by \eqref{expomix2} and $\gamma'$ is any number strictly greater than $\gamma$ given by \eqref{gamma}.
	
	\noindent \textbf{\textup{A)}} Assume $\tyk{D}$ is $(w,v)$-SDL for some strictly positive constants $w$ and $v$. Then for $u_n=r+\frac1v\log N_n$ we have
	\begin{align*}
	\lim_{n\to\infty} \mu\para{M_{I_n}\leq u_n}=e^{-w e^{-v r}}.
	\end{align*}
	\textbf{\textup{B)}} Assume $\tyk{D}$ is $(w_1,w_2,v)$-DL for some strictly positive constants $w_1,w_2$ and $v$. Then for $u_n=r+\frac1v\log N_n$ we have
	\begin{align*}
	e^{-w_1e^{-v r}}\leq \varliminf_{n\to\infty} \mu&\para{M_{I_n}\leq u_n}\leq\varlimsup_{n\to\infty} \mu\para{M_{I_n}\leq u_n}\leq e^{-w_2e^{-v r}}.
	\end{align*}
\end{thm}

Set $V(z)=\tub{x:\tyk{D}(x)\leq z}$. As mentioned in the introduction, the starting point of the proof is to write
\begin{gather}
\begin{split}
\label{eq:maxintegral}
\mu\para{M_{I_n}\leq u_n}=\int \een_{V(u_n)}(a_{m_{\alpha_n}}x)\een_{V(u_n)}(a_{m_{\alpha_n+1}}x)\cdots \een_{V(u_n)}(a_{m_{\beta_n}}x).
\end{split}
\end{gather}
Before we can apply exponential mixing to the integral we need to appropriately approximate the characteristic functions with smooth functions.

\subsection{Smooth approximations of characteristic functions}
The smooth approximations will be constructed by convolving the characteristic function with a smooth function specifically chosen for the purpose. Let $\mu_G$ denote the Haar measure on $G$ from which $\mu$ on $X$ is induced. Recall that for functions $\phi:G\to\R$ and $\psi:X\to\R$, the convolution $\phi*\psi:X\to\R$ is defined by
\begin{align*}
(\phi*\psi)(x)=\int_G \phi(g)\psi(g^{-1}x)\, d \mu_G(g).
\end{align*}
The smooth function, to be denoted $\phi_{\epsilon}\in C^{\infty}(G)$, will be chosen as follows. Let $\epsilon>0$. Pick a coordinate chart $\sigma:\R^d\to G$ such that $\R^d\supset B(0,\epsilon)\subset\sigma^{-1}(B(e,\epsilon))$ where $e\in G$ denotes the identity. On $\R^d$ we pick a function $\phi\in C^{\infty}(\R^d)$ such that $\supp(\phi)\subset B(0,1)$. Define then $\hat{\phi}_{\epsilon}:\R^d\to\R$ by 
\begin{align*}
\hat{\phi}_{\epsilon}(y)=\frac{\phi(\epsilon^{-1}y)}{\int_{B(0,\epsilon)} \phi(\epsilon^{-1}y)\rho(y)\,d\lambda(y)},
\end{align*}
where $\rho$ denotes the density of $\mu_G$ with respect to the Lebesgue measure $\lambda$ on $\R^d$.
Finally define
\begin{align}
\label{phiepsi}
\phi_{\epsilon}:=\hat{\phi}_{\epsilon}\circ\sigma^{-1}.
\end{align}
This function is seen to satisfy that $\int_{G} \phi_{\epsilon}\, d\mu_G=1$ and $\supp(\phi_{\epsilon})\subset B\para{e,\epsilon}$.
Notice that by doing a change of variables we get
\begin{align*}
\int_{B(0,\epsilon)} \phi(\epsilon^{-1}y)\rho(y)\,d\lambda(y)=
\epsilon^{d}\int_{B(0,1)} \phi(y)\rho(\epsilon y)\,d\lambda(y)=O(\epsilon^d),
\end{align*}
as $\epsilon\to 0$. Hence 
\begin{equation}
\label{phihatasymp}
\hat{\phi}_{\epsilon}(x)= O(\epsilon^{-d}\phi(\epsilon^{-1}x))\enskip \text{ as }\enskip \epsilon\to 0.
\end{equation}

We will need the following easy lemma concerning the Sobolev norm of such convolutions.
\begin{lemma}[]
\label{prelim1}
Let $A\subset X$ be measurable.
\begin{enumerate}[label=(\roman*)]
	\setlength\itemsep{0.4em}
	\item\label{prelim1i} $\int \phi_{\epsilon} * \een_{A} = \mu(A).$
	\item\label{prelim1ii} For any $\zeta\in\skr{g}$: $D_{\zeta}(\phi_{\epsilon}*\een_A)=(D_{\zeta}\phi_{\epsilon})*\een_A$.
	\item\label{prelim1iii} $S_k(\phi_{\epsilon}*\een_A)\leq S_k(\phi_{\epsilon})$.%$\sqrt{\mu(A)}$
	\item\label{prelimiv} $S_k(\phi_{\epsilon})=O\para{\epsilon^{-\para{\frac{d}{2}+k}}}$ as $\epsilon\to 0$.
\end{enumerate}
\begin{proof}
\ref{prelim1i}: This is an easy consequence of Fubini's Theorem, the definition of convolution and the fact that the action of $G$ on $X$ is measure preserving. 

\ref{prelim1ii}: This follows by using the definition of the derivative along with a change of variables.

\ref{prelim1iii}: This follows by realizing that \ref{prelim1ii} leads to
\begin{align*}
D_{\zeta_1}\dots D_{\zeta_d} (\phi*\een_{A})=\para{D_{\zeta_1}\dots D_{\zeta_d}\phi} *\een_{A}.
\end{align*}
Taking the $L^2$-norm and using the Young inequality for convolutions gives
\begin{align*}
%\label{eq:DiffConvo}
\norm{\para{D_{\zeta_1}\dots D_{\zeta_d}\phi} *\een_{A}}_2&\leq \norm{D_{\zeta_1}\dots D_{\zeta_d}\phi}_1\sqrt{\mu(A)}\\
&\leq \norm{D_{\zeta_1}\dots D_{\zeta_d}\phi}_2.
\end{align*}
The estimate now follows from the definition of the Sobolev norm.

\ref{prelimiv}: This is proven by applying the chain rule to $\phi_{\epsilon}=\hat{\phi}_{\epsilon}\circ\sigma^{-1}$ and using \eqref{phihatasymp} to obtain an estimate of the form
\begin{equation*}
	\norm{D_{\zeta}\phi_{\epsilon}}_2\leq C\norm{\epsilon^{-(d+1)}\psi_{\zeta}(\epsilon^{-1}y)}_2
\end{equation*}
for some constant $C>0$ and some smooth $L_2$ function $\psi_{\zeta}:\R^d\to\R$. A change of variables gives the estimate for a derivative of degree one and the general case follows similarly.
\end{proof}
\end{lemma}

With our mixing estimate \eqref{expomix2} in mind it is natural to impose some regularity assumption on the Sobolev norm of the approximating functions. For $k\in \N$ and $C>0$ we say that $\psi\in \tyk{S}_k(X)$ is $(C,k)$-regular if
\begin{align*}
S_k(\psi)\leq C\sqrt{\norm{\psi}_1}.
\end{align*}
The following lemma is a slight variation of \cite{KleinMarg2} Theorem 1.1 (see also \cite{KleinMarg} Lemma 4.2) and the proof is identical. The key point is that we can find sequences of smooth functions $g_{n,\epsilon}, h_{n,\epsilon}$ that approximate $\een_{V(u_n)}$ and are all $(C,k)$-regular, where the constants $C$ and $k$ are independent of $n$.
\begin{lemma}[]
	\label{lem:[Approx]}
	For $\epsilon>0$, define the sets
	\begin{align*}
	V'(u_n,\epsilon)&=\tub{x\in V(u_n):\dd(x,\partial V(u_n))\geq \epsilon}\\
	V''(u_n,\epsilon)&=\tub{x\in X:\dd(x,V(u_n))\leq \epsilon}
	\end{align*}
	and the functions
	\begin{align*}
	g_{n,\epsilon}=\phi_{\epsilon}*\een_{V'\para{u_n,\epsilon}}\quad \text{and}\quad h_{n,\epsilon}=\phi_{\epsilon}*\een_{V''\para{u_n,\epsilon}}. 
	\end{align*}
	Then for any $k\in\N$ and any $\delta>0$ there exist $\epsilon(\delta)>0$ such that $g_{n,\epsilon(\delta)}$ and $h_{n,\epsilon(\delta)}$ are  $\para{C\epsilon^{-\para{\frac{d}{2}+k}},k}$-regular for some $C>0$ and satisfy 
	\begin{align*}
	g_{n,\epsilon(\delta)}\leq \een_{V(u_n)}\leq h_{n,\epsilon(\delta)}\leq 1
	\end{align*}
	and
	\begin{align*}
	\mu\para{V(u_n-\delta)}\leq \int g_{n,\epsilon(\delta)} \leq \int h_{n,\epsilon(\delta)} \leq \mu\para{V(u_n+\delta)}.
	\end{align*}
\end{lemma}
Going forward we will write $\epsilon:=\epsilon(\delta)$ to simplify notation.
For any integers $i_1< i_2$, let 
\begin{align*}
G_{(i_1,i_2)}(x)=\prod_{i=i_1}^{i_2}g_{n,\epsilon}(a_{m_{i}}x)\quad\text{and}\quad H_{(i_1,i_2)}(x)=\prod_{i=i_1}^{i_2}h_{n,\epsilon}(a_{m_{i}}x).
\end{align*}
Then Lemma \ref{lem:[Approx]} implies that
\begin{align}\label{smoothbound}
\int G_{(\alpha_n,\beta_n)} \leq \mu\para{M_{I_n}\leq u_{n}}\leq\int H_{(\alpha_n,\beta_n)}.
\end{align}
Now we can estimate $\mu\para{M_{I_n}\leq u_{n}}$ by applying exponential mixing to the two integrals bounding it. The application of exponential mixing is demonstrated in the next lemma.
\begin{lemma}[]
	\label{ExpoEstim}
	For any $\epsilon>0$ we get
	\begin{align}
	\label{MixDiffEstG}
	\num{\int G_{(\alpha_n,\beta_n)}-\para{\int g_{n,\epsilon}}^{N_n}}\ll \epsilon^{-\para{\frac{d}{2}+k}}\sum_{c=\alpha_n+1}^{\beta_n} e^{-\delta m_c}S_k\para{G_{(\alpha_n,c-1)}}
	\end{align}
	and
	\begin{align}
	\label{MixDiffEstH}
	\num{\int H_{(\alpha_n,\beta_n)}-\para{\int h_{n,\epsilon}}^{N_n}}\ll \epsilon^{-\para{\frac{d}{2}+k}} \sum_{c=\alpha_n+1}^{\beta_n} e^{-\delta m_c}S_k\para{H_{(\alpha_n,c-1)}}.
	\end{align}
	\begin{proof}
		
		The two estimates are proven identically, so we do the calculation only for the first.
		Set $\psi=g_{n,\epsilon}$. Then
		\begin{align*}
		G_{(\alpha_n,\beta_n)}=G_{(\alpha_n,\beta_n-1)}(x)\psi(a_{m_{\beta_n}}x).
		\end{align*} 
		By applying exponential mixing to this product we get the estimate
		\begin{align}
		\label{eq:1}
		\num{\int G_{(\alpha_n,\beta_n)}-\int G_{(\alpha_n,\beta_n-1)}\int \psi} \ll e^{-\delta m_{\beta_n}} S_k\para{G_{(\alpha_n,\beta_n-1)}}S_k(\psi).
		\end{align}
		We repeat the procedure by writing $G_{(\alpha_n,\beta_n-1)}(x)=G_{(\alpha_n,\beta_n-2)}(x)\psi(a_{m_{\beta_{n}-1}}x)$ and applying exponential mixing to get
		\begin{align*}
		\num{\int G_{(\alpha_n,\beta_n-1)}-\int G_{(\alpha_n,\beta_n-2)}\int \psi} \ll e^{-\delta m_{\beta_{n}-1}} S_k\para{G_{(\alpha_n,\beta_n-2)}}S_k(\psi).
		\end{align*}
		Inserting this in \eqref{eq:1} gives
		\begin{align*}
		\num{\int G_{(\alpha_n,\beta_n)}-\int G_{(\alpha_n,\beta_n-2)}\para{\int\psi}^2} &\ll e^{-\delta m_{\beta_n}} S_k\para{G_{(\alpha_n,\beta_n-1)}}S_k(\psi)\\
		&+ e^{-\delta m_{\beta_{n}-1}} S_k\para{G_{(\alpha_n,\beta_n-2)}}S_k(\psi)\int\psi.
		\end{align*}
		We continue rewriting $G_{(\alpha_n,\beta_n-i)}$ and applying exponential mixing for all $1\leq i\leq \beta_n-\alpha_n-1$. 
		Once the process terminates for $i=\beta_n-\alpha_n-1$ we insert all the estimates, one after the other, into \eqref{eq:1}. This gives
		\begin{align*}
		\num{\int G_{(\alpha_n,\beta_n)}-\para{\int\psi}^{N_n}} \ll S_k(\psi)\sum_{c=\alpha_n+1}^{\beta_n} e^{-\delta m_c}S_k\para{G_{(\alpha_n,c-1)}}\para{\int\psi}^{\beta_n-c}.
		\end{align*}
		We use $S_k(\psi)=O\para{\epsilon^{-\para{\frac{d}{2}+k}}}$ and $\int\psi\leq 1$ by Lemma \ref{prelim1} and \ref{lem:[Approx]} to get the final estimate
		\begin{align*}
		\num{\int G_{(\alpha_n,\beta_n)}-\para{\int\psi}^{N_n}} \ll \epsilon^{-\para{\frac{d}{2}+k}}\sum_{c=\alpha_n+1}^{\beta_n} e^{-\delta m_c}S_k\para{G_{(\alpha_n,c-1)}}.	
		\end{align*}
		The implicit constant depends on $k$. However, $k$ remains fixed throughout and hence the dependence does not play a role here. We leave it out in the notation to reflect this. This will be the case in subsequent estimates as well.
 		\end{proof}
\end{lemma}
Having split the integral of a product into a product of integrals we are left with estimating the error term, i.e. the right hand side of \eqref{MixDiffEstG} and \eqref{MixDiffEstH}. % that are explicit in $n$. 
For this we need to understand the Sobolev norms of the functions $G_{(\alpha_n,c-1)}$ and $H_{(\alpha_n,c-1)}$ for $\alpha_n+1\leq c\leq \beta_n$. This is the content of Lemma \ref{SoboProd}. First we need the following estimate. 
\begin{lemma}
	\label{AppLem}
	Let $\psi$ be either of $g_{n,\epsilon}$ or $h_{n,\epsilon}$ and let $\zeta\in B$. For any $\epsilon>0$, any $q\in\N$ and any $\gamma'>\gamma$ we have
	\begin{align}
	\label{1stDiffEsti}
	\norm{D^q_{\Ad(a_{m_j})\zeta}\psi}_2\ll \epsilon^{-\para{\frac{d}{2}+q}} e^{q\gamma' m_j}.
	\end{align}
	for $j$ sufficiently large, where $\gamma$ is given by \eqref{gamma}.
\end{lemma}
To prove this lemma it is sufficient to consider the integer part of $m_j$ for which one can use that $B$ is a Jordan basis for $\Ad(a_1)$. Considering powers of the Jordan matrix and using linearity of the differential operator delivers the estimate. Since the proof consists of standard computations which are lengthy but well known, we leave out the details. From here on $\gamma'$ will generally denote any number strictly greater than $\gamma$.

\begin{lemma}[]
	\label{SoboProd}
	For any $\epsilon>0$ we have
	\begin{align*}
	S_k\para{G_{(\alpha_n,c-1)}}\ll \epsilon^{-k\para{\frac{d}{2}+1}} (c-\alpha_n)^k e^{k\gamma' m_{c-1}}
	\end{align*}
	and
	\begin{align*}
	S_k\para{H_{(\alpha_n,c-1)}}\ll \epsilon^{-k\para{\frac{d}{2}+1}} (c-\alpha_n)^k e^{k\gamma' m_{c-1}}
	\end{align*}
	for any $\alpha_n+1\leq c\leq\beta_n$.
	\begin{proof}
		Again, both estimates are proven identically and we give the calculation only for the first. Also, we give details of the calculation only for the first degree Sobolev norm, that is, $S_1(G_{(\alpha_n,c-1)})$. The argument for $S_k(G_{(\alpha_n,c-1)})$ follows almost identically and we explain the generalization at the end of the proof. 
		
		Set $\psi:=g_{n,\epsilon}$, $\psi_j(x):=\psi(a_{m_j}x)$ and recall that $\psi\leq 1$. Set also 
		\begin{align*}
		G:=G_{(\alpha_n,c-1)}=\psi_{\alpha_n}\psi_{\alpha_n+1}\cdots\psi_{c-1}.
		\end{align*}
		To obtain the upper bound we need an upper bound on $\Vert D_{\zeta_r}G\Vert_2$, $\zeta_r\in B$ which does not depend on the choice of $1\leq r\leq d$. For a given $\zeta_r\in B$ we first rewrite $D_{\zeta_r}G$. The product rule gives  
		\begin{gather}
		\begin{split}
		\label{soboprod1}
		D_{\zeta_r} G(x)&
		=\sum_{s=\alpha_n}^{c-1} \psi_{\alpha_n}(x)\cdots \frac{d}{dt}\kpara{\psi_{s}(\exp(t\zeta_r)x)}|_{t=0} \cdots  \psi_{c-1}(x).
		\end{split}
		\end{gather}
		We can rewrite $\psi_{s}(\exp(t\zeta_r)x)$ as follows.
		\begin{align*}
		\psi_{s}(\exp(t\zeta_r)x)
		&=\psi (a_{m_s}\exp(t\zeta_r)x)\\
		&=\psi(a_{m_s}\exp(t\zeta_r)a_{-m_s}a_{m_s}x)\\
		&=\psi(\exp(t\Ad(a_{m_s})\zeta_r)a_{m_s}x).
		\end{align*}
		Hence
		\begin{align*}
		\frac{d}{dt}\kpara{\psi_{s}(\exp(t\zeta_r)x)}|_{t=0}=D_{\Ad(a_{m_s})\zeta_r}\psi(a_{m_s} x).
		\end{align*}
		Inserting this in \eqref{soboprod1} and taking the $L^2$-norm then gives
		\begin{gather}
		\begin{split}
		\label{L2Esti}
		\Vert D_{\zeta_i} G\Vert_2&\leq\sum_{s=\alpha_n}^{c-1} \norm{\psi_{\alpha_n}(x)\cdots \frac{d}{dt}\kpara{\psi_{s}(\exp(t\zeta_i)x)}|_{t=0} \cdots  \psi_{c-1}(x)}_2\\
		&\leq \sum_{s=\alpha_n}^{c-1} \norm{\psi}_{\infty}^{c-\alpha_n}\norm{D_{\Ad(a_{m_s})\zeta_i}\psi}_2\\
		&\ll \epsilon^{-\para{\frac{d}{2}+1}} (c-\alpha_n)e^{\gamma'{m_{c-1}}},
		\end{split}
		\end{gather}
		where we used Lemma \ref{AppLem}, that $m_s\leq m_{c-1}$ and that $\psi\leq 1$. Finally,
		\begin{align*}
		S_1(G)^2=\Vert G\Vert_2^2+\sum_{n_1=1}^{d}\Vert D_{\zeta_{n_1}} G\Vert_2^2&\ll \Vert G\Vert_2^2+\sum_{n_1=1}^{d}\para{\epsilon^{-\para{\frac{d}{2}+1}}\para{c-\alpha_n} e^{\gamma' m_{c-1}}}^2.
		\end{align*}
		It follows that
		\begin{align*}
		S_1(G)\ll \epsilon^{-\para{\frac{d}{2}+1}} (c-\alpha_n)e^{\gamma' m_{c-1}}.
		\end{align*}
		The estimate for $S_k(G)$ runs almost identical. In this case we need an upper bound on $\norm{D_W G}_2 =\Vert D_{\zeta_1}^{n_1}\cdots D_{\zeta_d}^{n_d} G\Vert_2$ which does not depend on the choice of $n_1,\dots,n_d$, but only on the sum $n_1+\cdots+n_d=\deg(D_W)$. Set $l:=\deg(D_W)$. From the product rule we know that $D_W G$ is a sum consisting of $(c-\alpha_n)^l$ terms. Each term is a product of $(c-\alpha_n)$ functions. Some of these functions will be of the form $\psi_i$ while others will be derivatives of $\psi_i$ of some degree between 1 and $l$ where $\alpha_n\leq i\leq c-1$ . The important common property of each term is that the total degree of derivatives in the term is $l$. When estimating the $L_2$ norm of such a product we first bound each factor by its own $L_{\infty}$ norm to obtain a product of $L_{\infty}$ norms. Analogue to \eqref{L2Esti} we bound the $L_{\infty}$ norms of the not differentiated $\psi_i$'s by 1. For the $L_{\infty}$ norms of the differentiated $\psi_i$'s, we apply the Young inequality for convolutions. This states for general functions that $\norm{f*g}_{\infty}\leq \norm{f}_2\norm{g}_2$. By the definition of the $\psi_i$'s and Lemma \ref{prelim1}\ref{prelim1ii} we see that
		\begin{equation*}
		 \norm{D^q_{\Ad(a_{m_s})\zeta}\psi_i}_{\infty}\leq \norm{D^q_{\Ad(a_{m_s})\zeta} \phi_\epsilon}_2=O\para{\epsilon^{-\para{\frac{d}{2}+q}} e^{q\gamma' m_s}}
		\end{equation*}
		by Lemma \ref{AppLem}.
		This upper bound is maximized when the product contains $l$ different $\psi_i$'s all differentiated one time. This gives the estimate    
		\begin{align*}
		\norm{D_W G}_2 =\Vert D_{\zeta_1}^{n_1}\cdots D_{\zeta_d}^{n_d} G\Vert_2=O\para{\epsilon^{-l\para{\frac{d}{2}+1}}(c-\alpha_n)^l e^{l\gamma' m_{c-1}}}
		\end{align*}
		and by the definition of the Sobolev norm we get the final estimate
		\begin{align*}
		S_k(G)\ll \epsilon^{-k\para{\frac{d}{2}+1}} (c-\alpha_n)^k e^{k\gamma' m_{c-1}}.
		\end{align*}
	\end{proof}
\end{lemma}
\begin{remark}
	One might be tempted to simplify this proof by using an estimate of the form $S_k(f_1\cdots f_m)\ll_{k} S_{k+r}(f_1)\cdots S_{k+r}(f_m)$ for some $r>0$ (see for example \cite{BEG} (1.13)). However, for large $m$ this generates a worse estimate, hence the additional work is justified.
\end{remark}

Combining Lemma \ref{ExpoEstim} and Lemma \ref{SoboProd} we get the following estimates.
\begin{align}
\label{newerrorg}
\num{\int G_{(\alpha_n,\beta_n)}-\para{\int g_{n,\epsilon}}^{N_n}}\ll C_{\epsilon} \sum_{c=\alpha_n+1}^{\beta_n} e^{k\gamma' m_{c-1}-\delta m_c}(c-\alpha_n)^k
\end{align}
and
\begin{align}
\label{newerrorh}
\num{\int H_{(\alpha_n,\beta_n)}-\para{\int h_{n,\epsilon}}^{N_n}}\ll C_\epsilon \sum_{c=\alpha_n+1}^{\beta_n} e^{k\gamma' m_{c-1}-\delta m_c}(c-\alpha_n)^k.
\end{align}
with $C_\epsilon:=\epsilon^{-\frac{d}{2}(k+1)-2k}$. This completes the preparations for the proof of Theorem \ref{thm:[GenSDL]}.

\subsection*{Proof of Theorem \ref{thm:[GenSDL]}}
\label{Larg}

In the following we prove part \textbf{A)} of the theorem. Part \textbf{B)} is almost identical and we make a comment on this at the end of the proof. 

Assume $\tyk{D}$ to be $(w,v)$-SDL for some $w,v>0$. 
By applying \eqref{newerrorg} and \eqref{newerrorh} to \eqref{smoothbound} we get,
\begin{align}
\label{newstyleapprox}
\para{\int g_{n,\epsilon}}^{N_n}-O_{\epsilon}&\para{\sum_{c=\alpha_n+1}^{\beta_n} e^{k\gamma' m_{c-1}-\delta m_c}(c-\alpha_n)^k}\nonumber\\
\leq \mu\bigg (&M_{I_n}\leq u_{n}\bigg)\\
&\leq \para{\int h_{n,\epsilon}}^{N_n}+O_{\epsilon}\para{\sum_{c=\alpha_n+1}^{\beta_n} e^{k\gamma' m_{c-1}-\delta m_c}(c-\alpha_n)^k}.\nonumber
\end{align}
The constant $C_{\epsilon}$ will not play a role in this proof hence we integrated it in the implicit constant. We now determine the limit of the upper and lower bound. We start with the error term. 
The assumption that $\sup_{j\in\N}\frac{m_{j-1}}{m_j}<\min\para{1,\frac{\delta}{k\gamma'}}$ implies the existence of a $\rho>1$ such that for all $j\in\N$
\begin{align*}
\frac{m_{j-1}}{m_j}\leq\frac{1}{\rho} <\min\para{1,\frac{\delta}{k\gamma'}}.
\end{align*}
it follows that
\begin{align*}
m_j\geq \rho m_{j-1}\geq \dots \geq \rho^{j-1}m_{1}.
\end{align*}
Here we stopped at $\rho^{j-1}m_{1}$ instead of $\rho^j m_0$ since $m_1>0$ while $m_0$ could be zero. It follows that
\begin{align*}
k\gamma' m_{j-1}-\delta m_j=m_j\para{k\gamma'\frac{m_{j-1}}{m_j}-\delta}\leq m_{1}\rho^{-1} \para{\frac{k\gamma'}{\rho}-\delta}\rho^{j}<0.
\end{align*}
Set $\sigma=-m_{1}\rho^{-1}\para{\frac{k\gamma'}{\rho}-\delta}>0$. Then, since $\alpha_n\to\infty$ for $n\to\infty$, we see that for sufficiently large $n\in\N$ we have
\begin{align*}
\sum_{c=\alpha_n+1}^{\beta_n} e^{k\gamma' m_{c-1}-\delta m_c}(c-\alpha_n)^k&\leq\sum_{c=\alpha_n+1}^{\beta_n} \para{e^{-\sigma}}^{\rho^c}(c-\alpha_n)^k\\
&\leq \sum_{c=\alpha_n+1}^{\infty} \para{e^{-\sigma}}^{\rho^c}(c-\alpha_n)^k.
\end{align*}
This series is convergent and hence the right hand side goes to 0 as $n\to\infty$. This proves that
\begin{align}
\label{errorvanish}
\sum_{c=\alpha_n+1}^{\beta_n} e^{k\gamma' m_{c-1}-\delta m_c}(c-\alpha_n)^k\to 0 \quad\text{for}\quad n\to\infty.
\end{align}
Returning to the main terms, we claim that
\begin{align*}
%\label{eq:easyclaim}
\lim_{\delta\to 0}\lim_{n\to\infty}\para{\int g_{n,\epsilon}}^{N_n}=\lim_{\delta\to 0}\lim_{n\to\infty}\para{\int h_{n,\epsilon}}^{N_n}=e^{-we^{-v r}},
\end{align*}
To see this we first look at the limits for $n\to\infty$. Recall from Lemma \ref{lem:[Approx]} that for any $\delta>0$ we can find $\epsilon=\epsilon(\delta)>0$ such that
\begin{align}
\label{eq:subsetsAgain}
\mu\para{V(u_n-\delta)}\leq \int g_{n,\epsilon} \leq \int h_{n,\epsilon} \leq \mu\para{V(u_n+\delta)},
\end{align}
and hence
\begin{align*}
\mu\para{V(u_n-\delta)}^{N_n}\leq \para{\int g_{n,\epsilon}}^{N_n} \leq \para{\int h_{n,\epsilon}}^{N_n} \leq \mu\para{V(u_n+\delta)}^{N_n}.
\end{align*}
Recall that $u_n=r+\frac1v\log N_n$. Using the $(w,v)$-SDL property of $\tyk{D}$ we get
\begin{align*}
\mu\para{V(u_n\pm\delta)}^{N_n}&=\para{1-\frac{w e^{-v(r\pm\delta)}}{N_n}+o\para{\frac{e^{-v(r\pm\delta)}}{N_n}}}^{N_n}\\
&=e^{N_n \log\para{1-\frac{w e^{-v(r\pm\delta)}}{N_n}+o\para{\frac{e^{-v(r\pm\delta)}}{N_n}}}}.
\end{align*}
We estimate the right hand side using the second order Taylor expansion $\log (1+x)=x+O(x^2)$. Inserting the resulting estimate and taking the limit for $n\to\infty$ gives
\begin{align*}
e^{w e^{-v(r-\delta)}}\leq \lim_{n\to\infty} \para{\int g_{n,\epsilon}}^{N_n}\leq \lim_{n\to\infty}\para{\int h_{n,\epsilon}}^{N_n}\leq e^{w e^{-v(r+\delta)}}.
\end{align*}
Therefore the claim is proved by letting $\delta\to 0$. 
In conclusion we have proved that \textbf{A)} follows by taking limits $n\to\infty$ and then $\delta\to 0$ in \eqref{newstyleapprox}.

The proof of part \textbf{B)} is essentially identical to the proof of part \textbf{A)}. The only difference occurs following equation \eqref{eq:subsetsAgain} where we apply the DL property instead of the SDL property. For some strictly positive constants $w_1, w_2$ and $v$ we then get
\begin{align*}
1-\frac{w_1 e^{-v(r-\delta)}}{N_n}&\leq \int g_{n,\epsilon}\leq \int h_{n,\epsilon} \leq 1-\frac{w_2 e^{-v(r+\delta)}}{N_n},
\end{align*}
which gives the inequalities claimed in \textbf{B)} when raising to the power $N_n$ and taking the $\varliminf$ and $\varlimsup$ for $n\to\infty$ followed by $\delta\to 0$.

\section{Closest distance returns}
\label{ClosRet}

In this section we prove Theorem \ref{riedistret}. Assume throughout this section that 
\begin{equation}
\label{obsd}
\tyk{D}(x)=-\log\dd(x,x_0)
\end{equation}
for some $x_0\in X$. The following theorem is a restatement of Theorem \ref{riedistret} with the constant $C$ made explicit.
\begin{thm}[]
	\label{thm:[ClosRet]}
	Let $\rho=\para{\sup_{s\in\N}\frac{m_{s-1}}{m_s}}^{-1}$. Assume that $m_j$ satisfies 
	\begin{align*}
	\sup_{j\in\N}\frac{m_{j-1}}{m_j}<\min\para{1,\frac{\delta}{k\gamma'}}
	\end{align*}
	where $k,\delta$ are as in \eqref{expomix2} and $\gamma'$ is any constant strictly bigger than $\gamma$ as defined in \eqref{gamma}. Assume further that
	\begin{align}
	\label{techcond}
	N_n=o\para{e^{\sigma \rho^{\alpha_n}}},
	\end{align}
	where $\sigma=\frac{-m_{1}\rho^{-1}}{k\para{\frac32+\frac{2}{d}}+\frac12}\para{\frac{k\gamma'}{\rho}-\delta}>0$. Then there exists a constant $\kappa>0$ such that for $u_n(r)=r+\frac1d\log N_n$ and all $x_0\in X$ we have
	\begin{align*}
	\lim_{n\to\infty} \mu\para{M_{I_n}\leq u_n(r)}=e^{-\kappa e^{-d r}}.
	\end{align*}
\end{thm}
The proof follows the same strategy as the proof of Theorem \ref{thm:[GenSDL]}. Asymptotic estimates for the tail distribution function of $\tyk{D}$ are easy in this case since $\mu\para{x:\tyk{D}(x)\geq z}=\mu(B(x_0,e^{-z}))$ and it is well known that
\begin{align}
\label{ballasymp}
	\mu(B(x_0,e^{-z}))=\kappa e^{-dz}+o(e^{-dz})\quad\text{ as } z\to\infty.
\end{align}
for some constant $\kappa>0$ since $\dd$ is a Riemannian metric and $\mu$ is the Haar measure.

As previously stated, it is the fact that $\tyk{D}$ is not uniformly continuous that forces some changes. In particular, we must adapt Lemma \ref{lem:[Approx]}, i.e. the way we approximate characteristic functions by smooth functions. Recall the notation $V(r)=\tub{x:\tyk{D}(x)\leq r}$ and note that $\een_{V(u_n)}=\een_{B(x_0,e^{-u_n})^c}$. In the following lemma we approximate $\een_{V(u_n)}$ by $(C_n,k)$-regular functions, i.e. we are forced to allow $C$ to depend on $n$.
\begin{lemma}[]
	\label{retapprox}
	For $\omega>0$, set $\epsilon_{n,\omega}:=e^{-u_n}(1-e^{-\omega})$ and let $\phi_{\epsilon_{n,\omega}}$ be defined as in \eqref{phiepsi}. Define the functions
	\begin{align*}
	g_{n,\omega}=\phi_{\epsilon_{n,\omega}}*\een_{B(x_0,e^{-u_n+ \omega})^c}\quad \text{and}\quad h_{n,\omega}=\phi_{\epsilon_{n,\omega}}*\een_{B(x_0,e^{-u_n- \omega})^c}.
	\end{align*}
	For any $k\in \N$, there exists $C>0$ such that for all $\omega>0$ and $n\in\N$, $g_{n,\omega}$ and $h_{n,\omega}$ are  $\para{C\epsilon_{n,\omega}^{-(\frac{d}{2}+k)},k}$-regular functions such that 
	\begin{align*}
	g_{n,\omega}\leq \een_{V(u_n)} \leq h_{n,\omega}\leq 1
	\end{align*}
	and
	\begin{align}
	\label{retapproxeq2}
	\mu\para{B(x_0,e^{-u_n+ \omega})^c}= \int g_{n,\omega}\leq \int h_{n,\omega}=\mu\para{B(x_0,e^{-u_n- \omega})^c}.
	\end{align}
	\begin{proof}
		Write $\epsilon:=\epsilon_{n,\omega}$. For any $\omega>0$ we have
		\begin{align*}
		B(x_0,e^{-u_n+ \omega})^c\subset B(x_0,e^{-u_n})^c \subset B(x_0,e^{-u_n - \omega})^c.
		\end{align*}
		Then \eqref{retapproxeq2} follows directly from Lemma \ref{prelim1}\ref{prelim1i}. We also get
		\begin{align}
		\label{eq:conc22}
		g_{n,\omega}\leq \een_{B(x_0,e^{-u_n})^c}\leq h_{n,\omega}\leq 1.
		\end{align}
		To see this, we first rewrite $g_{n,\omega}$ as an integral, 
		\begin{align*}
		g_{n,\omega}(x)&=\int_G \phi_{\epsilon}(g)\een_{B(x_0,e^{-u_n+ \omega})^c}(g^{-1}x)\,d m(g)\\
		&=\int_{B(e,\epsilon)} \phi_{\epsilon}(g)\een_{gB(x_0,e^{-u_n+ \omega})^c}(x)\,d m(g).
		\end{align*}
		Its clear that if $gB(x_0,e^{-u_n+ \omega})^c\subset B(x_0,e^{-u_n})^c$ then the first inequality of \eqref{eq:conc22} is established. So assume $x\in gB(x_0,e^{-u_n+ \omega})^c$. Then we can write $x=gy$ where $y\in B(x_0,e^{-u_n+ \omega})^c$. This means that $\dd(y,\partial B(x_0,e^{-u_n})^c)\geq e^{-u_n+\omega}-e^{-u_n}=e^{-u_n}(e^{\omega}-1)\geq \epsilon$. Hence $gy\in B(x_0,e^{-u_n})^c$ and $gB(x_0,e^{-u_n+ \omega})^c\subset B(x_0,e^{-u_n})^c$. The second inequality of \eqref{eq:conc22} is proved similarly and the third inequality is trivial.
		
		The $\para{C\epsilon_{n,\omega}^{-(\frac{d}{2}+k)},k}$-regularity of $g_{n,\omega}$ and $h_{n,\omega}$ follows directly from Lemma \ref{prelim1} \ref{prelim1iii} and \ref{prelimiv}. 
	\end{proof}
\end{lemma}

\subsubsection*{Proof of Theorem \ref{thm:[ClosRet]}}
\vspace{0.4cm}
The proof is a copy of the proof of Theorem \ref{thm:[GenSDL]} with small changes. We obtain an analog of \eqref{newstyleapprox} which is
\begin{align}
\label{newstyleapprox2}
\para{\int g_{n,\omega}}^{N_n}&-O\para{C_{\epsilon_{n,\omega}}\sum_{c=\alpha_n+1}^{\beta_n} e^{k\gamma' m_{c-1}-\delta m_c}(c-\alpha_n)^k}\nonumber\\
\leq \mu\bigg (M_{I_n}&\leq u_{n}\bigg)\\
\leq &\para{\int h_{n,\omega}}^{N_n}+O\para{C_{\epsilon_{n,\omega}}\sum_{c=\alpha_n+1}^{\beta_n} e^{k\gamma' m_{c-1}-\delta m_c}(c-\alpha_n)^k}.\nonumber
\end{align}
Again, we want to determine the limit of the upper and lower bound as $n\to\infty$ and $\omega\to 0$ and we begin by looking at the error term. First we make the trivial observation 
\begin{align*}
C_{\epsilon_{n,\omega}}\sum_{c=\alpha_n+1}^{\beta_n} e^{k\gamma' m_{c-1}-\delta m_c}(c-\alpha_n)^k \leq C_{\epsilon_{n,\omega}}(\beta_n-\alpha_n)^k\sum_{c=\alpha_n+1}^{\beta_n} e^{k\gamma' m_{c-1}-\delta m_c}.
\end{align*}
We now look at the sum. Since $\sup_{j\in\N}\frac{m_{j-1}}{m_j}=\frac{1}{\rho}$ we see that for all $j\in\N$
\begin{align*}
m_j\geq \rho m_{j-1}\geq \dots \geq \rho^{j-1}m_{1},
\end{align*}
and consequently,
\begin{align*}
k\gamma' m_{j-1}-\delta m_j=m_j\para{k\gamma'\frac{m_{j-1}}{m_j}-\delta}\leq m_{1}\rho^{-1}\para{\frac{k\gamma'}{\rho}-\delta}\rho^{j}<0.
\end{align*}
Set $\lambda=-m_{1}\rho^{-1}\para{\frac{k\gamma'}{\rho}-\delta}>0$. Now, since $\rho>1$, we know that for $j\in\N$ sufficiently large we have that $\rho^{j+i}\geq \rho^j + i$ and hence
\begin{align*}
\para{e^{-\lambda}}^{\rho^{j+i}}\leq \para{e^{-\lambda}}^{\rho^{j}+i}.
\end{align*}
Consequently, for $n\in\N$ sufficiently large we can write
\begin{gather}
\begin{split}
\label{sumconv2}
\sum_{c=\alpha_n+1}^{\beta_n} e^{k\gamma' m_{c-1}-\delta m_c}\leq\sum_{c=\alpha_n+1}^{\beta_n} \para{e^{-\lambda}}^{\rho^c}&\leq \sum_{i=\rho^{(\alpha_n+1)}}^{\rho^{(\alpha_n+1)}+(\beta_n-\alpha_n)} \para{e^{-\lambda}}^{i}\\
&\leq \sum_{i=\rho^{\alpha_n}}^{\infty} \para{e^{-\lambda}}^{i}\\
&=O\para{e^{-\lambda\rho^{\alpha_n}}}.
\end{split}
\end{gather}
We now look at the rest of the error term. Recalling from Lemma \ref{retapprox} that $\epsilon_{n,\omega}=e^{-u_n}(1-e^{-\omega})$, we get
\begin{gather}
\begin{split}
\label{epsilonn}
\epsilon_{n,\omega}^{-\frac{d}{2}(k+1)-2k}(\beta_n-\alpha_n)^k&\leq\para{e^{-r}e^{-\frac{1}{d}\log N_n}(1-e^{-\omega})}^{-\frac{d}{2}(k+1)-2k}N_n^k\\
&= O_{r,\omega}\para{N_n^{k\para{\frac32+\frac{2}{d}}+\frac12}}.
\end{split}
\end{gather}
Hence for $n\in\N$ sufficiently large we have,
\begin{align}
\label{eq:reterrorvan}
 \epsilon_{n,\omega}^{-\frac{d}{2}(k+1)-2k}(\beta_n-\alpha_n)^k\sum_{s=\alpha_n+1}^{\beta_n} e^{k\gamma' m_{s-1}-\delta m_s}
=O\para{\frac{N_n^{k\para{\frac32+\frac{2}{d}}+\frac12}}{e^{\lambda\rho^{\alpha_n}}}}\to 0
\end{align}
for $n\to\infty$ since we assumed that $N_n=o\para{e^{\sigma \rho^{\alpha_n}}}$ for
\begin{equation*}
	\sigma=\frac{-m_{1}\rho^{-1}}{k\para{\frac32+\frac{2}{d}}+\frac12}\para{\frac{k\gamma'}{\rho}-\delta}=\frac{\lambda}{k\para{\frac32+\frac{2}{d}}+\frac12}.
\end{equation*}
For the main terms we claim that
\begin{align}
\label{eq:easyclaim2}
\lim_{\omega\to 0}\lim_{n\to\infty}\para{\int g_{n,\omega}}^{N_n}=\lim_{\omega\to 0}\lim_{n\to\infty}\para{\int h_{n,\omega}}^{N_n}=e^{-\kappa e^{-d r}}.
\end{align}
The proof is again almost identical to the corresponding part of the proof of Theorem \ref{thm:[GenSDL]}. It follows from Lemma \ref{retapprox} and \eqref{ballasymp}, that
\begin{align*}
\int g_{n,\omega}=\mu\para{B(x_0,e^{-(u_n-\omega)})^c}=1-\kappa\para{e^{-d(u_n-\omega)}}+o\para{e^{-d(u_n-\omega)}}\\
\int h_{n,\omega}=\mu\para{B(x_0,e^{-(u_n+\omega)})^c}=1-\kappa\para{e^{-d(u_n+\omega)}}+o\para{e^{-d(u_n+\omega)}}.
\end{align*}
Again, insert $u_n(r)=r+\frac1d\log N_n$, raise to the power $N_n$, apply the second order Taylor expansion of $\log(1+x)$ and take limits for $n\to\infty$ to obtain 
\begin{align*}
e^{-\kappa e^{d\omega}e^{-d r}}= \lim_{n\to\infty} \para{\int g_{n,\omega}}^{N_n}\leq \lim_{n\to\infty}\para{\int h_{n,\omega}}^{N_n}= e^{-\kappa e^{-d\omega}e^{-d r}}.
\end{align*}
The claim is then proved by letting $\omega\to 0$. It then follows from \eqref{eq:reterrorvan} and \eqref{eq:easyclaim2} that when we take the limit for $n\to\infty$ and then for $\omega\to 0$ in \eqref{newstyleapprox2}, we get
\begin{align*}
\lim_{n\to\infty}\mu\para{M_{I_n}\leq u_{n}}
= e^{-\kappa e^{-d r}}.
\end{align*}

\end{document}